\documentclass[%
 reprint,
 amsmath,amssymb,
 aps,
]{revtex4-2}

\usepackage{graphicx}
\usepackage{dcolumn}
\usepackage{xcolor}
\usepackage{bm}
\usepackage{subfigure}

\usepackage{amssymb}
\usepackage{enumitem}
\usepackage{algorithm2e}
\usepackage{adjustbox}

\newcommand{\vect}[1]{\boldsymbol{#1}}

\begin{document}

\preprint{APS/123-QED}


\title{Machine Learning for the identification of phase-transitions \\ in interacting agent-based systems: a Desai-Zwanzig example}

\author{Nikolaos Evangelou$^{1,2}$, Dimitris G. Giovanis$^{3, 4}$,  \\ George A. Kevrekidis$^2$, Grigorios A. Pavliotis$^{5}$, Ioannis G. Kevrekidis$^{1,2,\star}$}

\affiliation{$^1$Department of Chemical and Biomolecular Engineering, Johns Hopkins University, 3400 North Charles Street, Baltimore, MD 21218, USA}
\affiliation{$^2$Department of Applied Mathematics and Statistics, Johns Hopkins University, 3400 North Charles Street, Baltimore, MD 21218, USA}
\affiliation{$^3$Department of Civil and Systems Engineering, Johns Hopkins University, 3400 North Charles Street, Baltimore, MD, 21218, USA}
\affiliation{$^4$Hopkins Extreme Materials Institute, Johns Hopkins University,  3400 North Charles Street, Baltimore, MD 21218, USA}
\affiliation{$^5$Department of Mathematics, Imperial College London, London SW7 2AZ, United Kingdom}

\email{yannisk@jhu.edu}




\begin{abstract}
Deriving closed-form, analytical expressions for reduced-order models, and judiciously choosing the closures leading to them, has long been the strategy of choice for studying phase- and noise-induced transitions for agent-based models (ABMs). 
In this paper, we propose a data-driven framework that pinpoints phase transitions for an ABM- the Desai-Zwanzig model - in its mean-field limit, using a smaller number of variables than traditional closed-form models. 
To this end, we use the manifold learning algorithm Diffusion Maps to identify a parsimonious set of data-driven latent variables, and show that they are in one-to-one correspondence with the expected theoretical order parameter of the ABM. 
We then utilize a deep learning framework to obtain a conformal reparametrization of the data-driven  coordinates that facilitates, in our example, the identification of a single parameter-dependent ODE in these coordinates. We identify this ODE through a residual neural network inspired by a numerical integration scheme (forward Euler). We then use the identified ODE -- enabled through an odd symmetry transformation-- to construct the bifurcation diagram exhibiting the phase transition.
\end{abstract}

\keywords{interacting particle systems, machine-learning, dynamical systems, phase-transitions, system identification}                
\maketitle


\section{\label{sec:Introduction}Introduction:}

Complex dynamic phenomena are ubiquitous in natural sciences, social sciences, and engineering \cite{naldi2010mathematical,pareschi2013interacting,helfmann2021statistical,Helfmann2021}.
In many cases, their study has been performed using agent-based models (ABM) also known as interacting particle systems (IPS). 
The dynamics of those models in the thermodynamic (mean-field) limit undergo phase transitions -- bifurcations in non-linear dynamics terminology ~\cite{ChayesPanferov2010, CGPS2020}. 
The exploration of the dynamics of these systems typically involves extensive numerical simulation scenarios for a very large number, $N$, of interacting agents that can be truly challenging and impedes the widespread utilization of these models.
%
Therefore, coarse-graining methodologies become necessary in order to reduce this complexity. Reducing the dimensionality of an ABM can be achieved by defining collective variables that are capable of accurately describing the full dynamics of such large systems in terms of a relatively small number of observables \cite{gross2006epidemic,winkelmann2021mathematical,zagli2023dimension}; or discovering those by data mining \cite{liu2014coarse,liu2015equation,fabiani2023tasks}. 

In our previous work \cite{zagli2023dimension}, we developed and successfully tested a model reduction approach based on the cumulants of the single-agent probability distribution for such large ABM systems. 
%
The proposed coarse-graining framework was based on an analytical closure methodology of the infinite hierarchy of equations for the moments or, equivalently, cumulants of the probability distribution of the infinite-dimensional system. The basic steps for building a reduced-order Desai-Zwanzig (DZ) model are: (i) Consider the mean-field ansatz by writing the $N$-particle distribution function, the solution of the $N-$particle Fokker-Planck equation, as the product of one-particle distribution functions $\rho(x,t)$  \cite{MartzelAslangul2001}. 
(ii) Represent $\rho(x,t)$ either in terms of its moments (DZ~\cite{Dawson1983}) or in terms of its Fourier coefficients (Smoluchowski/noisy Kuramoto model~\cite{RevModPhys.77.137}), thus obtaining an infinite system of ordinary differential equations (ODEs) that is {\emph exactly} equivalent to the McKean-Vlasov PDE. Truncating the equations for the moments (or the Fourier coefficients), while choosing an appropriate closure scheme gives the reduced-order model described and analyzed in \cite{zagli2023dimension}.  However, selecting the correct closure scheme for the truncated moments' system is the \textit{trickiest} part of the coarse-graining procedure. 

\begin{figure*}[ht!]
    \centering
 \includegraphics[width = 15cm]{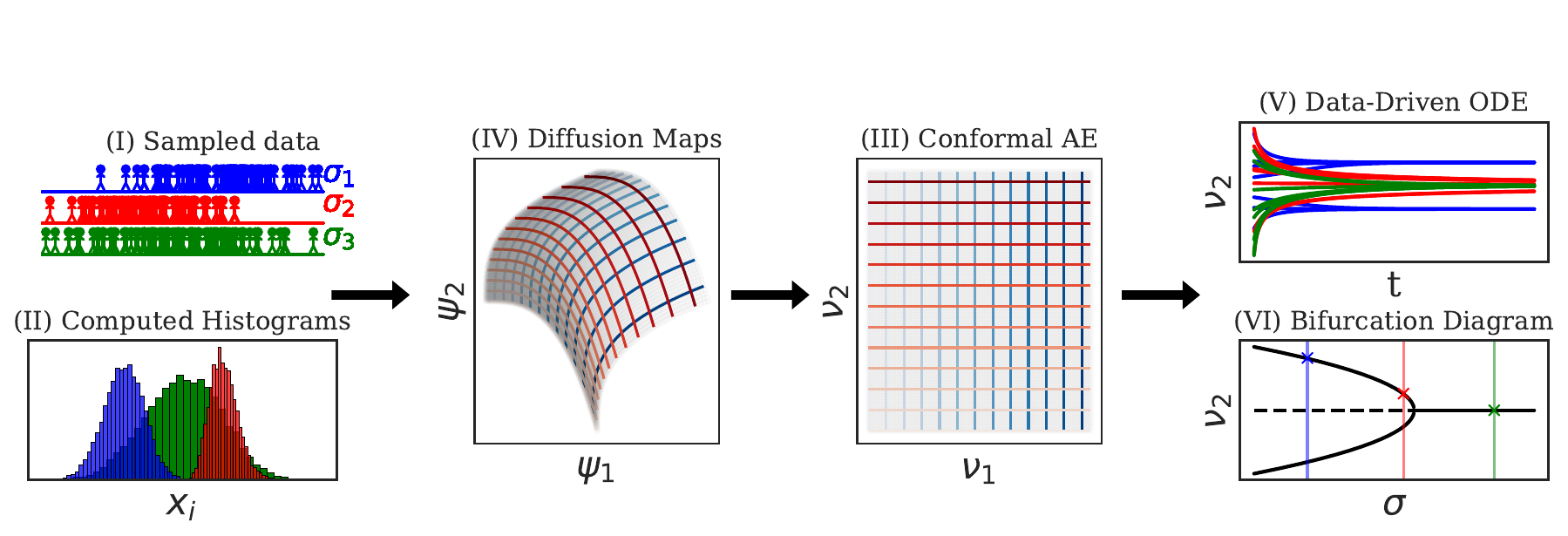}
    \caption{Schematic of the overall workflow. (I) Sample data from the Agent-Based Model (ABM) across multiple initial conditions and parameter values, (II) Compute histograms for each snapshot of the ABM. (III) Apply the Diffusion Maps algorithm on the computed histograms to discover a reduced latent embedding. (IV) Use a conformal autoencoder (AE) to find a conformal reparametrization of the latent space. (V) Identify a data-driven ODE in terms of the latent coordinate $\nu_2$ of the AE. (VI) Construct the bifurcation diagram (enabled via a symmetry transformation) in terms of the latent coordinate $\nu_2$.}
    \label{fig:intro-figure}
\end{figure*}

%
This analytical approximation requires choosing the ``correct'' observable(s), the right
level of observation/analytical closure, and --importantly -- the level at which we attempt the closure.  To overcome these constraints, this paper proposes a data-driven framework for studying phase transitions of ABMs by (i) discovering the coarse observables in a data-driven fashion, (ii) determining the level of closure, and (iii) identifying the reduced dynamics directly from data, (iv) utilizing the data-driven reduced model and applying an odd symmetry transformation to it for the construction of the bifurcation diagram.
To this end, we employ Diffusion Maps \cite{coifman2006diffusion}, a manifold learning algorithm, conformal autoencoders \cite{evangelou2022parameter} and residual neural networks inspired by numerical integrators of ordinary differential equations \cite{rico1992discrete,gonzalez1998identification,dietrich2021learning,evangelou2023learning,fabiani2023tasks}. This allows us to circumvent the difficulties arising from the choice of an analytical model and the selection of its closure.  The main steps of our proposed data-driven workflow, are illustrated in Figure \ref{fig:intro-figure}.

As an example, we apply our data-driven framework to the generalization of the DZ model with multiplicative noise considered in~\cite{zagli2023dimension}; this very well-studied ABM model features a second-order (continuous) phase transition, which enables us to compare the performance of our data-driven approach with well-known analytical and computational results.  
We emphasize, however, that our method is quite general, since it does not depend on the detailed features of the interaction, and can be applied to many different ABMs that exhibit phase transitions in the thermodynamic limit.

\section{\label{sec:Model}The Agent-Based Model:}
\noindent
We consider a model describing an ensemble of $N$ identical interacting agents subject to multiplicative noise. For this model, the agents are coupled via a mean reverting force and a second-order phase transition exists at a critical temperature 
that can be calculated analytically~\cite{zagli2023dimension}[Eqn. A12]. This transition appears as a symmetric pitchfork bifurcation of the mean-field dynamics. The dynamics of each agent, $x_i$, are described by a stochastic differential equation (SDE) of the form
\begin{equation}
\label{eq:ABM_final}
    \mathrm{d}x_i=\bigg[ -x_i^3 + (\alpha + \nu\sigma_m^2)x_i - \theta(x_i - \overline{x}) \bigg]dt + \sqrt{\sigma^2 + \sigma_m^2x_i^2} \mathrm{d}W_i;
\end{equation}
where $\sigma$ is the bifurcation parameter, $\theta$ denotes the interaction strength, $\alpha$ is a parameter characterizing the amplitude of the multiplicative noise, $\nu$ corresponds to different mathematical prescriptions of the SDEs (It\^o, Stratonovich, etc.) and, $\sigma_m$ is a rectifying parameter: when $\sigma_m \neq 0$ the phase-transition is pushed to higher values of $\sigma.$ The agents are coupled through $\bar{x}$, which denotes the center of mass of the system (equal to the first moment $M_1$). Furthermore, $dW_i, i=1 \dots N$ denotes independent Brownian motions. The values of the parameters were set to $\alpha = 1$,  
 $\theta = 4$ , $\sigma_m = 0.8$, $\nu = \frac{1}{2}$, and $N = 12,000$ as in \cite{zagli2023dimension}.
 \par 
For this model, it was demonstrated in \cite{zagli2023dimension} that, away from the phase transition at $\sigma \cong 1.890$, the first moment, i.e., the order parameter of the system, is sufficient for an accurate description of the mean-field dynamics.  As we get closer to the phase transition/the bifurcation, more moments need to be taken into account in order to accurately locate the bifurcation.
More specifically, it was numerically demonstrated that at least a four-moment truncation is necessary in order to accurately predict the bifurcation/phase transition. In this work, by observing the eigenvalues of the Jacobian at the steady state (based on the moments equation proposed in \cite{zagli2023dimension}), we confirm that a clear separation of time scales prevails in the neighborhood of the bifurcation, and that the long-term dynamics live on a one-dimensional manifold. Our goal here, discussed in the next section, is to find a data-driven parametrization of this one-dimensional manifold, and identify the dynamics on it. 




\section{\label{sec:Results}Numerical Results:}
\subsection{Latent data-driven coordinates}

\begin{figure*}[htpb!]
    \centering
    \subfigure[]{ \includegraphics[width=0.15\textwidth]{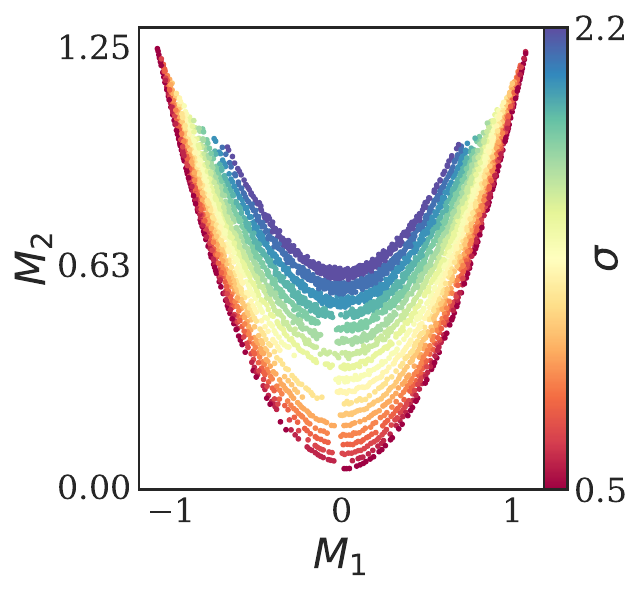}
    \label{fig:moments_colored}}
\subfigure[]{ \includegraphics[width=0.135\textwidth]{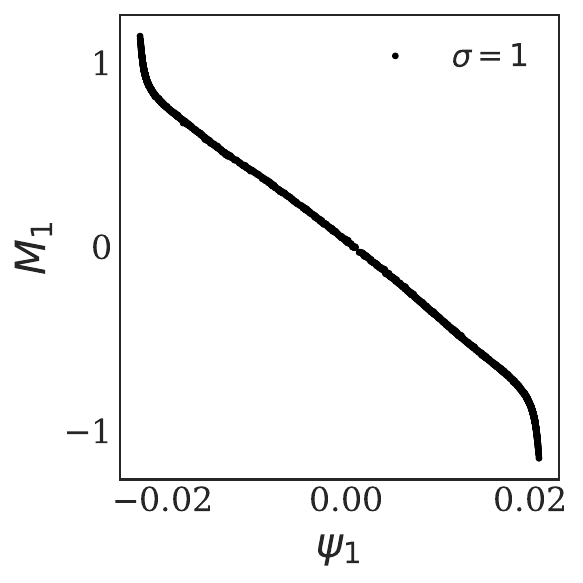}
\label{fig:one_to_one_psi_m1}}  
\subfigure[]{
{\includegraphics[clip,width=0.14\textwidth]{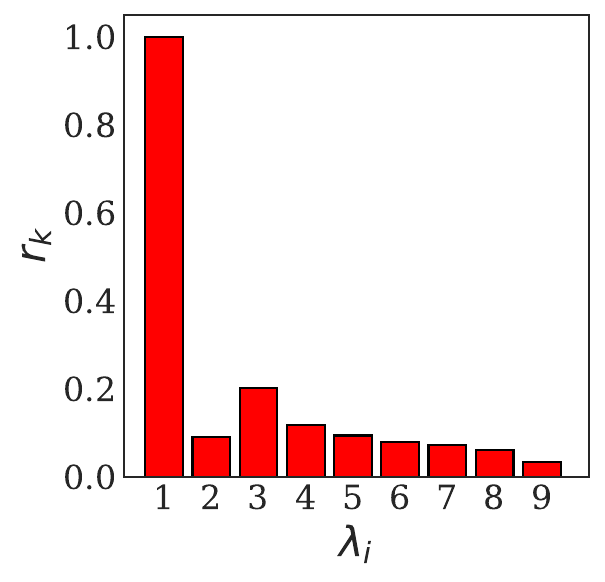}}\label{fig:local_linear_1_sigam}}
    \subfigure[]{ \includegraphics[width=0.15\textwidth]{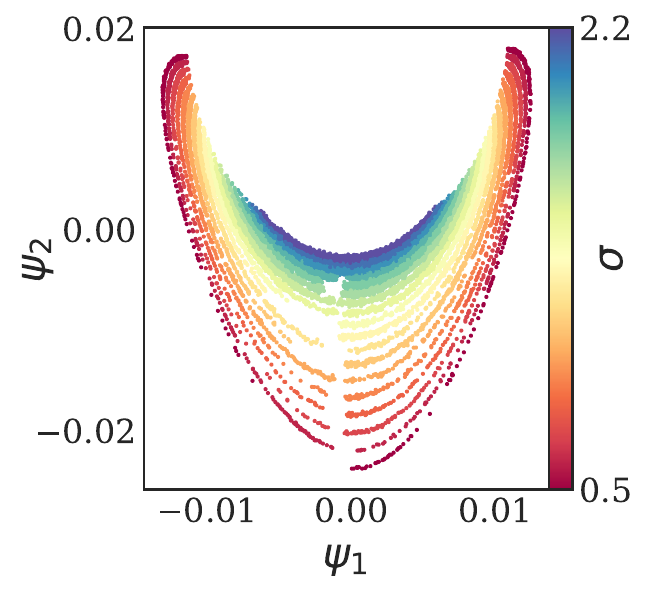}  \label{fig:non_harmonic_sigma_hist}}
 \subfigure[]{ 
\includegraphics[width=0.15\textwidth]{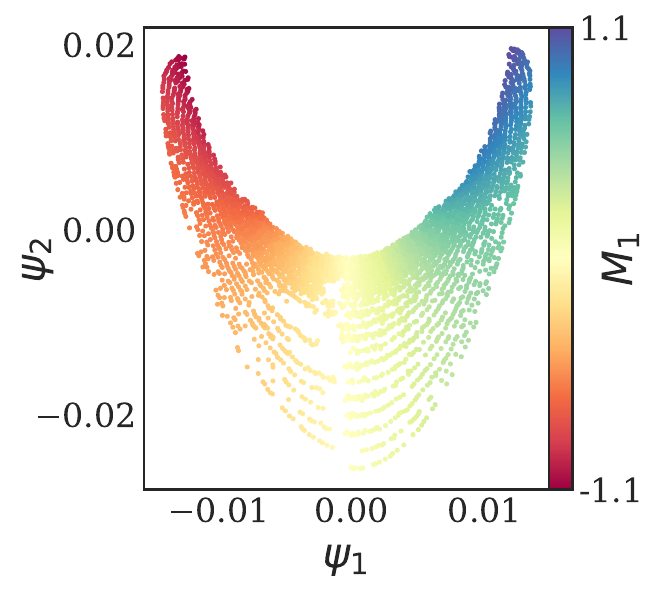}
\label{fig:non_harmonics_m1_hist}}
\subfigure[]{
\includegraphics[width=0.14\textwidth]{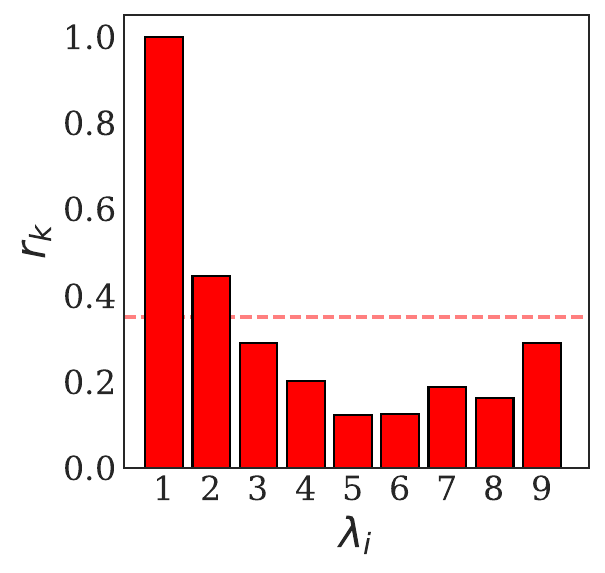}
\label{fig:local_linear_regression_hist}}
    \caption{
    (a) The first two moments ($M_1$, $M_2$) estimated from the ABM simulations colored with the parameter value $\sigma$ at which they were obtained. (b-c) Diffusion Maps on histograms for a single value of the parameter $\sigma$: (b) One-to-one relation between leading histogram moment $M_1$ and leading Diffusion Maps coordinate $\psi_1$; (c) shows the residual ($r_{k}$) based on the local linear regression algorithm indicating that $\psi_1$ might be enough to parametrize the data (larger value of $r_k$). (d-f) Diffusion Maps on collected histograms from the ABM simulation across multiple values of the parameter $\sigma$: (d-e) show the non-harmonic diffusion map coordinates 
    $(\psi_1, \psi_2)$ colored with $\sigma$ and $M_1$ respectively; (f) shows the residual ($r_{k}$) based on the local linear regression algorithm indicating that $\psi_1$ and $\psi_2$ might be enough to parametrize the data (larger values of $r_k$). 
    }
\label{fig:data_diffusion_maps}
\end{figure*}

We start by sampling data on an equidistant grid of $18$ distinct values of the parameter $\sigma$ in the interval $\sigma \in [0.5,2.2]$. A more detailed description of the sampling strategy is discussed in Section \ref{sec:data_collection}. 

Given the sampled data from the ABM we compute the first four moments $M_1,M_2,M_3,M_4$ as
    \hbox{$M_k = \frac{1}{N} \sum_{i=1}^N x_i^k(t)$},
where $x_i^k$ denotes the $i$-th agent at a fixed time raised in the power $k$, $M_k$ is the $k$-th moment of the population from the sampled ABM data. The computation of the moments serves as a comparison between our data-driven framework and the one proposed in \cite{zagli2023dimension} where an analytical system of ODEs based on the four moments is proposed. 
The choice of four moment is also supported by studying the time-scale separation of the analytical reduced model proposed in \cite{zagli2023dimension} 
 between the four-moment truncation and higher-order moment truncations of the dynamics, see \ref{sec:sep_time_scales} in the Appendix. In Figure \ref{fig:moments_colored} we plot the first two moments $M_1$ and $M_2$ colored with the parameter $\sigma$. Based on this, we argue that the data are at least two-dimensional, and that the parameter $\sigma$ appears to correlate with the moments $M_1,M_2$. 

We then proceed by computing Diffusion Maps (see Section \ref{sec:Dmaps} for a description of the algorithm) using data from the sampled ABM simulation. The Diffusion Maps algorithm was computed using: (a) histograms obtained \textit{for a single value} of the parameter $\sigma$, (b) histograms obtained across our grid of 18 distinct $\sigma$ values, (c) on the four computed moments from data sampled across our 18 values of $\sigma$; those computations are reported in Section  \ref{sec:dmaps_SI} of the Appendix. We would like to reiterate that the Diffusion Maps computations on the moments are not necessary for the overall framework they serve as a comparison of our approach and the one by \cite{zagli2023dimension} et al.


Computing Diffusion Maps (based on histograms) for single values of $\sigma$ (far from the phase transition) gave, a one-dimensional manifold parametrized by $\psi_1$. In Figure \ref{fig:one_to_one_psi_m1} we plot the leading Diffusion Maps coordinate $\psi_1$ against the first moment $M_1$ and demonstrate that they are one-to-one. This suggests that Diffusion Maps is capable of discovering a coordinate that is one-to-one with the known order parameter $M_1$ for this model~\cite{Dawson1983,PhysRevE.49.2639, zagli2023dimension}. The claim that the manifold is one-dimensional in this case is also supported from the results of the local-linear algorithm shown in Figure \ref{fig:local_linear_1_sigam}; the residual $r_k$ of the first eigenvector $\psi_1$ is larger than the one of the remaining. In Section \ref{sec:dmaps_SI} of the Appendix we illustrate the relation of the Diffusion Maps coordinate with the subsequent moments $M_2,M_3,M_4$.

We proceed with the Diffusion Maps computation on data sampled over multiple parameter values of $\sigma$. 
Details for the Diffusion Maps algorithm applied on histograms are provided in Section \ref{sec:dmaps_ABM}. Here in Figures \ref{fig:non_harmonic_sigma_hist}, \ref{fig:non_harmonics_m1_hist} we show the two, non-harmonic eigenvectors \textit{that parameterize different eigendirections} $\psi_1$  and $\psi_2$ colored with the parameter $\sigma$ and the first moment $M_1$. The selection of the non-harmonic eigenvectors was made by using the local linear regression algorithm proposed in \cite{dsilva2018parsimonious} and implemented by the \textit{datafold} Python package \cite{lehmberg2020datafold}. The residual $r_k$ of the first two eigenvectors, shown in Figure \ref{fig:local_linear_regression_hist}, 
 is greater than the residual of the remaining eigenvectors. This suggests that the first two eigenvectors ($\psi_1,\psi_2$) are non-harmonic, while the remaining ones are harmonic (functions) of the first two. 

The Diffusion Maps results for the case of using collected data across multiple parameter values, shown in Figures 
\ref{fig:local_linear_1_sigam},
\ref{fig:non_harmonic_sigma_hist}, and \ref{fig:local_linear_1_sigam}, suggest that the mean-field dynamics of the ABM live on a two-dimensional 
 manifold. However, the parameter $\sigma$ appears to correlate with the behavior of the moments observed at that $\sigma$ value. This suggests that even though the data can be parametrized by two coordinates $\psi_1$, $\psi_2$,  
 one might be able to \textit{disentangle} this coupled two-dimensional description into a factorized (a single state)$\times$(a parameter) description. In the next section, we illustrate how a data-driven reparametrization of the Diffusion Maps coordinates can be achieved that will allow us to \textit{disentangle} the effect of the parameter $\sigma$ from that of the latent variable.

\subsection{Y-shaped conformal autoencoder}

\begin{figure*}[ht!]
    \centering
    \begin{tabular}{cc}
    \subfigure[]{
    \raisebox{-0.45\height}{%
        \includegraphics[trim={1.75cm 1.75cm 1.75cm 1.75cm},clip,width=0.38\textwidth]{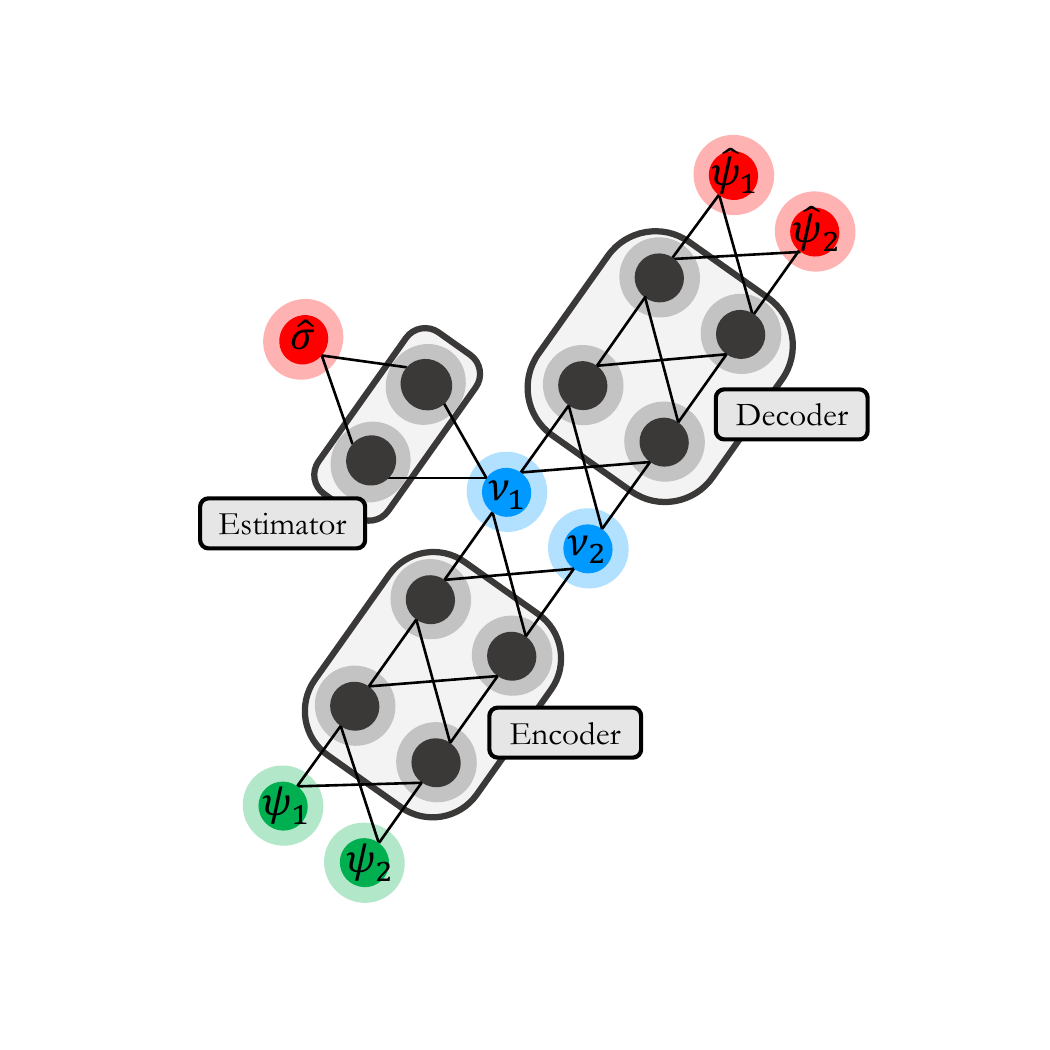}%
    }%
    \label{fig:y_shaped_nn}}
    \begin{tabular}{cc}
    \subfigure[]{
\includegraphics[width=0.19\textwidth]{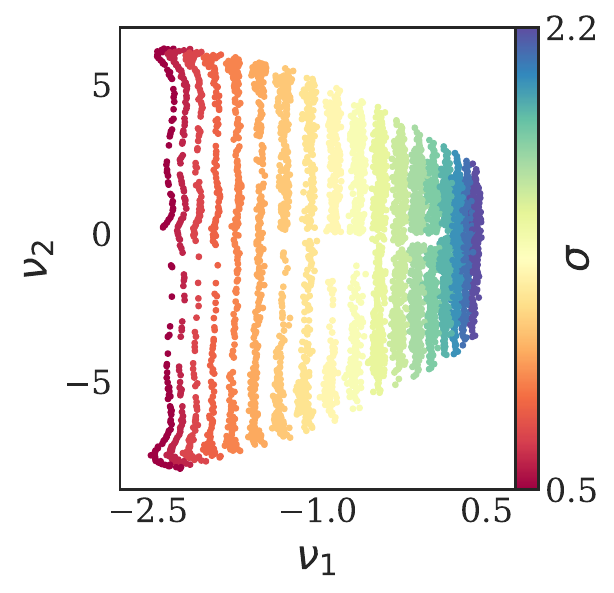}
\label{fig:cae_coordinates_sigma}}
    \subfigure[]{
\includegraphics[width=0.19
\textwidth]{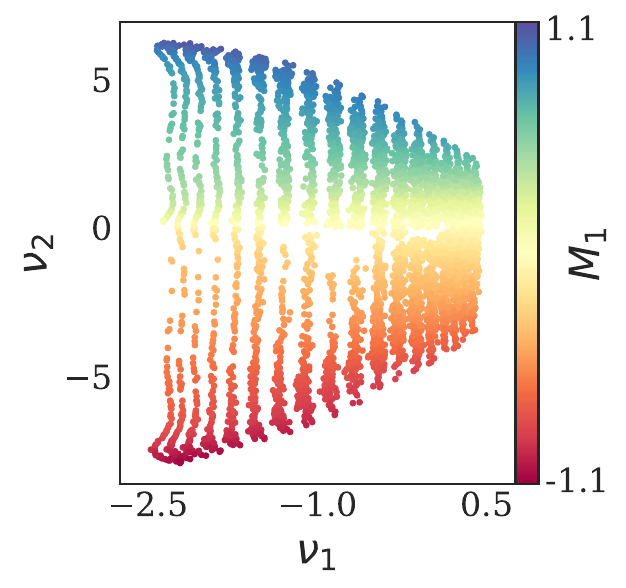}
\label{fig:cae_coordinates_m1}}
\\
    \subfigure[]{
\includegraphics[width=0.19\textwidth]{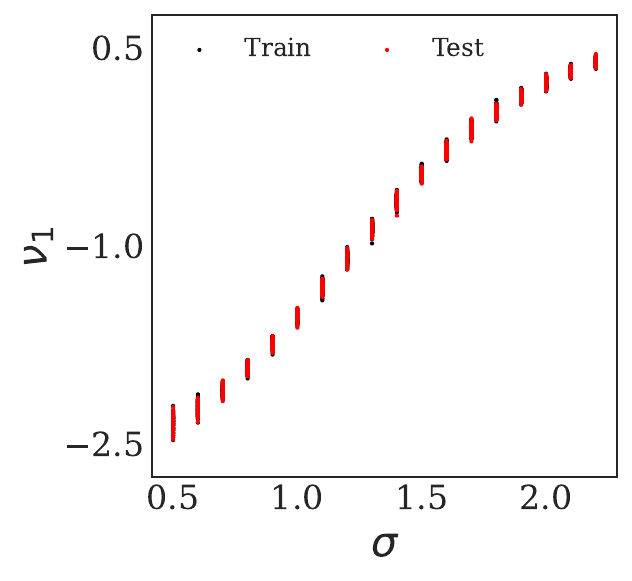}
    \label{fig:nu1_vs_sigma}}
    \subfigure[]{
\includegraphics[width=0.18\textwidth]{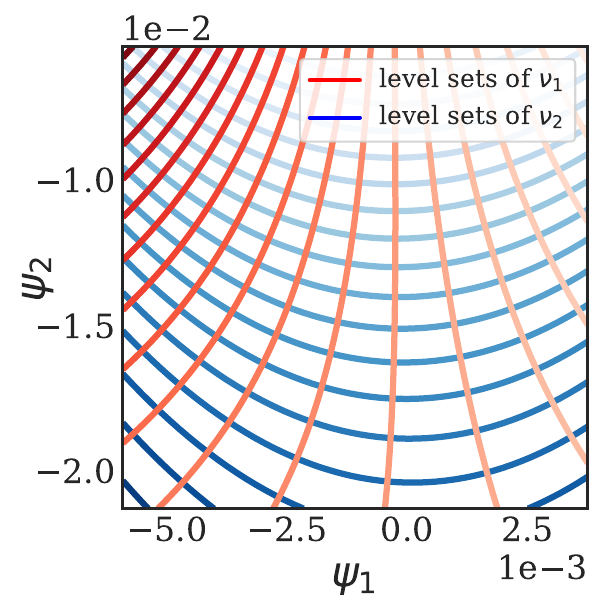}
    \label{fig:CAE_conformality}}
    \end{tabular}            
    \end{tabular}
    \caption{(a) A schematic of the Y-shaped conformal autoencoder. The inputs to the network ($\psi_1$ and $\psi_2$) are shown as green nodes. The outputs of the autoencoder ($\hat{\psi}_1$ and $\hat{\psi}_2$) and the estimator are shown as red nodes. The latent variables ($\nu_1$ and $\nu_2$) are shown as light blue nodes. (b-c) The obtained latent coordinates $\nu_1,\nu_2$ colored with $\sigma$ and $M_1$ respectively. We can see that the  $\sigma$ appears to vary only across $\nu_1$, compared to Fig.~\ref{fig:data_diffusion_maps} (c) in which $\sigma$ varied across both $\psi_1, \psi_2$. Similarly, $M_1$ appears to vary only across $\nu_2$. (d) The parameter $\sigma$ is plotted against the latent coordinate $\nu_1$ indicating a strong dependence. (e) Contour lines representing level sets of $\nu_1$ and $\nu_2$ are plotted in the Diffusion Maps ($\psi_1$,$\psi_2$) providing a visual illustration of the obtained conformality. }
    \label{fig:y-shaped}
\end{figure*}

In this section, we illustrate the use of the Y-shaped conformal autoencoder, originally introduced in our previous work \cite{evangelou2022parameter} as a means of disentangling the effect of different parameter combinations on model outputs in the context of parameter nonidentifiability. The Y-shaped conformal autoencoder seeks a reparametrization of the Diffusion Maps coordinates in which it disentangles the effect of the parameter $\sigma$ from that of the latent coordinates of the autoencoder. A more detailed description of the overall framework, the different loss function components of the network, and the training procedure we followed are discussed in Section \ref{sec:Y_shaped_CAE_details}. A schematic of the Y-shaped conformal autoencoder is shown in Figure \ref{fig:y_shaped_nn}. 
The coordinates $\psi_1$ and $\psi_2$, obtained by computing the Diffusion Maps algorithm on estimated histograms of the ABM, were the ones used as input to the Y-shaped conformal autoencoder.

 The bottleneck latent variables $\nu_1,\nu_2$ are shown in Figures \ref{fig:cae_coordinates_sigma} and \ref{fig:cae_coordinates_m1} colored with the parameter $\sigma$ and the first moment $M_1$ (that is the known coarse variable for this model \cite{zagli2023dimension}). We see that the parameter $\sigma$ appears to be varying only along $\nu_1$, compared to Figure \ref{fig:non_harmonic_sigma_hist} in which $\sigma$ varied across both $\psi_1, \psi_2$. Similarly, $M_1$ varies only along $\nu_2$. The latent variable $\nu_1$ shows a strong correlation with the parameter $\sigma$, as depicted by Figure \ref{fig:nu1_vs_sigma}, indicating that they are effectively one-to-one. 
The  conformality efficiency of the autoencoder becomes visually pronounced in Figure \ref{fig:CAE_conformality} where level sets of $\nu_1$ and $\nu_2$ are plotted in the Diffusion Maps coordinates.  Furthermore, the average value of the cosine of the angle, $\cos\theta =  \frac{ \nabla\nu_1 \cdot \nabla \nu_2 }{||\nabla \nu_1|| ||\nabla \nu_2 || }$, calculated over the test set, we found it to be $5 \times 10^{-4}$.

This reparametrization of the latent coordinates is crucial since, as we show in the next section, it allows us to identify a single-state dimension ODE depending on $\sigma$ (rather than a needlessly two-dimensional ODE).

In Section \ref{sec:Y_shaped_CAE_details} we provide additional results obtained for the Y-shaped autoencoder illustrating the reconstruction of its input data ($\psi_1, \psi_2$). Remarkably, in addition, the Y-shaped autoencoder allows us to estimate $\sigma$ from previously unseen histogram observation through $\nu_1$.

\subsection{Identifying parameter dependent ODEs}
\label{sec:Euler_NN_one}
\begin{figure*}[ht!]
    \centering
\subfigure[]{
\includegraphics[trim={0.375cm 0.5cm 0.6cm 0.6cm},clip,width=0.55\textwidth]{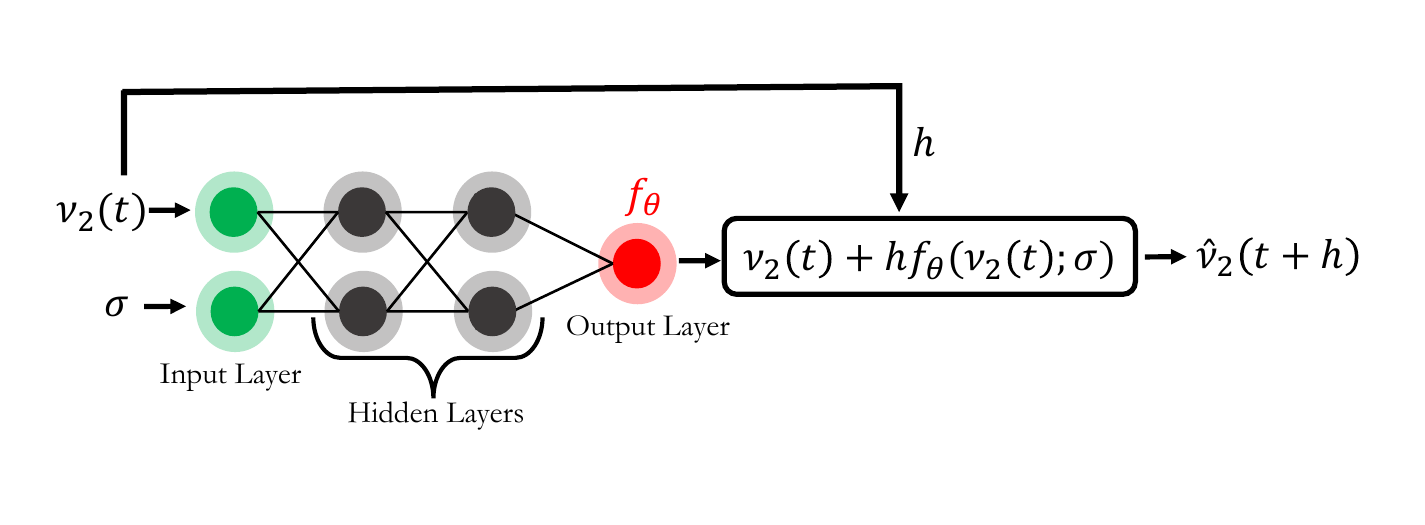}
\label{fig:euler_nn}
}
\subfigure[]{
\includegraphics[width = 0.18\textwidth]{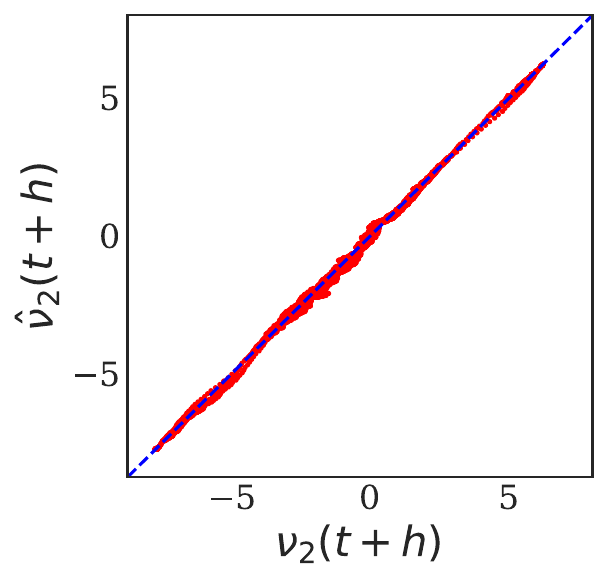}
\label{fig:7_parameters_parity_plot}}
\subfigure[]{
\includegraphics[width = 0.17\textwidth]{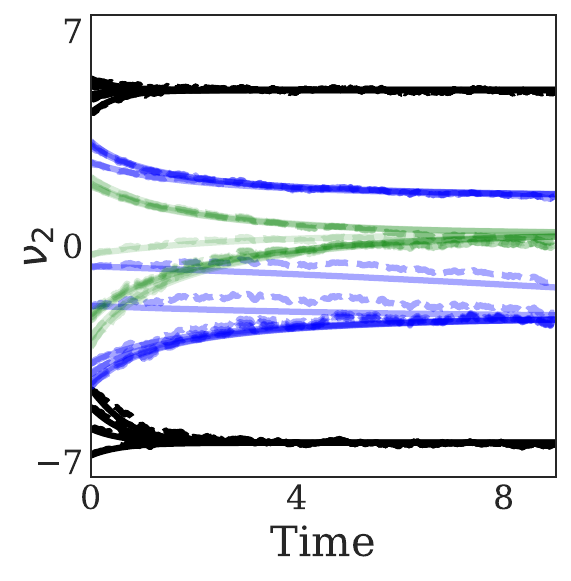}
\label{fig:3d_comparison}}
    \caption{(a) A schematic of the forward Euler residual neural network. The state variable $\nu_2(t)$, the parameter $\sigma$, and the time step $h$ are inputs to the neural network that estimates the right-hand side $f_{\theta}$ of the ODE. The right-hand-side is then used to estimate the state variable $\hat{\nu_2}(t + h)$ by using a forward Euler step.  (b) The predicted value for $\hat{\nu}_2(t +h)$ estimated from the Euler neural network is plotted against the true value (projected ABM trajectories in $\nu_2$ space) for values of the parameter $\sigma = \{0.57, 0.85, 1.11, 1.75, 1.9, 2.06, 2.25\}$ not included in the training set. (c) For three values, of the parameter not included in the training set we contrast generated paths from the ODE (solid lines) with paths simulated by the ABM (dashed lines) embedded in $\nu_2$. The colors black, blue (gray) and green (light-gray) corresponds to $\sigma = 1.11$, $\sigma = 1.75$ and $\sigma= 2.25$ respectively.}
\end{figure*}

Using the latent variable $\nu_2$ obtained from the Y-shaped conformal autoencoder as the state variable, we now identify a $\sigma$-dependent ODE; remember that $\nu_2$ is conformal to $\nu_1$ (and therefore to $\sigma$). The identification of the parameter-dependent ODE was achieved through a residual neural network inspired by the forward Euler numerical integration scheme. A schematic of the neural network is shown in Figure \ref{fig:euler_nn}. The inputs to this network are the state variable at time $t$, $\nu_2(t)$, and the parameter $\sigma$, and the output is the state variable evolved to time $t+h$, $\nu_2(t+h)$. We provide a more detailed description of the forward Euler network and the constructed loss function scheme in Section \ref{sec:euler_nn_methodology}.  

We test the ability of the ODE to produce accurate paths (trajectories), in $\nu_2$ space, for seven unseen values of the parameter $\sigma=\{0.57, 0.85, 1.11, 1.75, 1.9, 2.06, 2.25\}$, see Section \ref{sec:data_collection} for more details. In Figure \ref{fig:7_parameters_parity_plot} we illustrate the predicted values $\hat{\nu}_2(t+h)$ from the Euler neural network against the ground-truth (projected trajectories in $\nu_2$ generated by the ABM) for all seven unseen values of the parameter $\sigma$. The estimated MSE in this case was MSE = $6 \times 10^{-4}$.  
In addition, 
 in Figure \ref{fig:3d_comparison} we provide a visual comparison between simulated by the Euler Neural Network trajectories and the ABM model (projected in $\nu_2$) for the test parameter values $\sigma=\{1.11, 1.75, 1.9\}$. The paths generated by the ABM and those generated by the neural network-identified ODE show excellent visual agreement across different values of $\sigma$ and different initial conditions. This suggests that the identified right-hand-side provides a successful approximation of the ground-truth dynamics.

\subsection{Bifurcation Diagram: phase transition}

\begin{figure}[ht!]
    \centering
\subfigure[]{
\includegraphics[width = 0.22\textwidth]{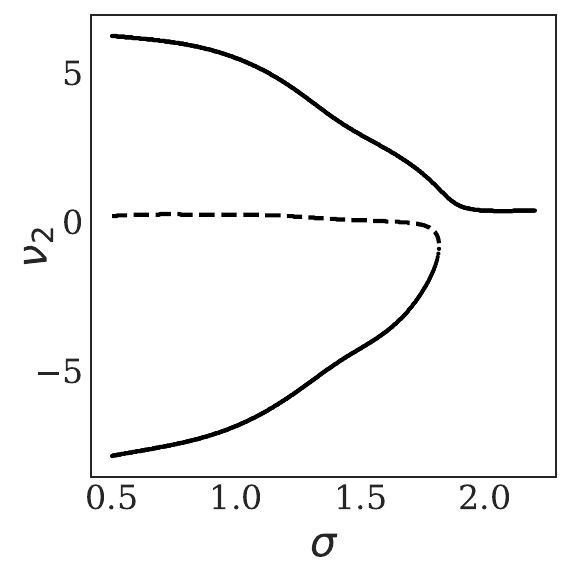}
\label{fig:bifurcation_disconnected}}
\subfigure[]{
\includegraphics[width = 0.22\textwidth]{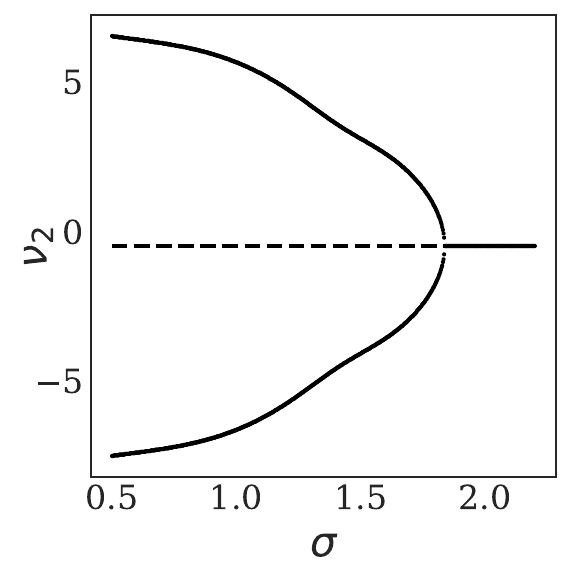}
\label{fig:bifurcation_connected}}
    \caption{
    (a) A representative bifurcation diagram constructed by the identified right-hand-side of the Euler neural network suggests a ``perturbed'' pitchfork. (b) A representative bifurcation diagram after applying equation \eqref{eq:symmetry} to the identified right-hand-side shows a symmetric pitchfork.}
    \label{fig:bifurcation_diagrams}
\end{figure}
In this section, we now proceed to test whether the parameter-dependent ODE identified in a data-driven way through our neural network can qualitatively and quantitatively capture the phase transition. To assess the robustness of our approach, we also estimate the error associated with the phase transition by training 5,000 neural networks in parallel, using different splits of the data for training, validation, and testing. All  hyperparameters were kept the same across models during this computation. For each of the trained networks given the identified right-hand-side, we compute the steady-states across a range of values of the parameter $\sigma$ and construct the bifurcation diagram. 
For values of $\sigma \in [0.5,2.2]$ a representative constructed bifurcation diagram is shown in Figure \ref{fig:bifurcation_disconnected}. The bifurcation diagram shows that the model identifies the existence of three steady states, two stable and one unstable, for $\sigma<1.84$. For $\sigma \geq 1.84$ a unique stable steady state exists. However, the computed bifurcation diagram clearly possesses the flip symmetry of the generic pitchfork bifurcation. We have identified a \textit{perturbed} pitchfork, in which the upper branch remains permanently stable and the lower branch exhibits a turning point where a stable and an unstable component collide. The consistency of identifying a perturbed pitchfork was maintained across the different models we trained.

To make sure we recover a generically symmetric pitchfork bifurcation diagram, a simple transformation suffices: given the identified right-hand-side $f_{\theta}(\nu_2(t);\sigma)$ of the neural network we compute
\begin{equation}
    \label{eq:symmetry}
    g(\nu_2;\sigma) = \frac{f_{\theta}(\nu_2(t);\sigma) - f_{\theta}(-\nu_2(t);\sigma)}{2}.
\end{equation}

\noindent
The vector field $g(\nu_2;\sigma)$ is then used to construct the bifurcation diagram, shown in Figure \ref{fig:bifurcation_connected} correctly capturing the symmetry.
Note that the computation of Equation \eqref{eq:symmetry} does not require a second neural network but only evaluations of the trained Euler neural network discussed in the previous section. 
The average critical parameter of the identified ODE (based on the $5000$ trained neural networks) was estimated at $\sigma^* = 1.837$ and the standard deviation was $0.026$; the true critical transition has been computed to be at $\sigma^* = 1.890$ \cite{zagli2023dimension}--see also\cite{Dawson1983}[Eqn. 3.52] for the case of additive noise. This suggests that the identified ODE provides a qualitatively correct and arguably quantitatively accurate approximation of the phase transition. This small discrepancy can likely be attributed to the limited range of $\sigma$ values used in our training data. With only $18$ distinct values for $\sigma$ (details on sampling in Section \ref{sec:data_collection}), the closest ones to the true bifurcation point are $\sigma = 1.80$ and $\sigma = 1.90$.  Using a finer grid of parameter values we believe it would have led to even smaller discrepancy.

\section{\label{sec:Conclusions}Summary and Conclusions}

We presented a data-driven framework for identifying qualitatively and approximating quantitatively phase transitions of interacting agent systems
through an interplay between direct simulations and machine-learning-assisted, coarse-grained system identification and bifurcation analysis.
We demonstrated that our framework is capable of identifying the phase transition of the DZ system with a single, interpretable data-driven variable, in contrast to the derived closed-form model, proposed in \cite{zagli2023dimension} where a system of four approximate ODEs were required for the task.
To discover the effective coarse (collective) variables that parameterize collected data from the ABM, the Diffusion Maps manifold learning algorithm \cite{coifman2006diffusion} was used. We showed that, for a constant parameter value ($\sigma =1$), away from the phase transition, the Diffusion Maps algorithm discovers a one-dimensional manifold parametrized by the data-driven coordinate $\psi_1$; and that $\psi_1$ is one-to-one with the theoretical order parameter $M_1$.
We then illustrated that, while data collected over a range of parameter values (bracketing the critical value) lie on a two-dimensional manifold,
this manifold can be ``disentangled" (factored) into a one-dimensional state-variability manifold ``cross'' a one-dimensional parameter variability manifold.
To disentangle the variability due to the parameter $\sigma$ from the state variability, in terms of the latent Diffusion Maps coordinates ($\psi_1,\psi_2$), we introduced a Y-shaped conformal autoencoder, initially proposed in \cite{evangelou2022parameter} for addressing parameter non-identifiability. We illustrated that the autoencoder's latent variables  ($\nu_1,\nu_2$) form a conformal set of coordinates that ultimately allow us to learn a single effective nonlinear ODE, in terms of a single state variable $\nu_2$ and depending on a single parameter $\sigma$ (the latter being one-to-one with $\nu_1$). This effective nonlinear ODE was identified via a residual neural network, templated on the forward Euler numerical integrator.
We compared generated paths of the identified ODE 
for test values of the parameter $\sigma$ to paths generated with the full ABM embedded in the latent coordinate $\nu_2$ and demonstrated good agreement between the two. To construct the bifurcation diagram from the identified right-hand-side of the ODE, we used an odd symmetry transformation of the right-hand-side. This transformation provides a symmetric vector field and captures the pitchfork bifurcation that denotes the phase transition. Imposing known physical symmetries in the data-driven models in order to enhance accuracy and improve generalizability is currently an active area of research \cite{olver2015modern,mattheakis2019physical,yarotsky2022universal,blum2022equivariant,villar2022dimensionless,olver2023normal, jin2020sympnets,alet2021noether,burby2020computing}.

Our proposed framework can be extended to a broad class of complex, multiscale dynamical systems for which a fine scale, atomistic/stochastic mode exists, but for which accurate closed form 
macroscopic equations at the coarse-grained level are not explicitly known. 
We are particularly interested in models from the social sciences, where no {\emph physics-informed} natural choice for the order parameter(s) exists. As an example, we mention the noisy Hegselmann-Krause model for opinion dynamics~\cite{Chazelle_al2017a}, for which it has been rigorously proved that a discontinuous disorder (no consensus)/order (consensus) phase transition exists~\cite{CGPS2020}[Prop. 6.2]. It would be interesting to identify the phenomenological order parameter introduced in~\cite{Chazelle_al2017a} using our approach. Other examples for which our approach is expected to be applicable are the Keller-Segel model for chemotaxis, see the work in \cite{lee2023learning,psarellis2022data,siettos2010system,setayeshgar2005application}, the Fitzhugh–Nagumo from Lattice Boltzmann in
\cite{armaou2005equation,galaris2022numerical}, the mean-field Fitzhugh–Nagumo model for Neurons  ~\cite{pavliotis2022method}[Sec. 5.6], and the Dean-Kawasaki stochastic PDE \cite{helfmann2021statistical,cornalba2023deankawasaki}, that takes into account finite-particle effects and fluctuations.


\begin{acknowledgments}
I.G.K. acknowledges partial support from the US AFOSR FA9550-21-0317 and the US Department
of Energy SA22-0052-S001. G.A.P. is partially supported by the Frontier Research Advanced Investigator Grant ERC grant Machine-aided general framework for fluctuating dynamic density functional theory. The authors are grateful to N. Zagli and A. Zanoni for useful discussions for making available the code from ~\cite{zagli2023dimension} and~\cite{pavliotis2022method}, respectively. Many thanks to T. Gaskin for helping set up the collaboration with Imperial College London. He acknowledges the hospitality of the Johns Hopkins group. 
\end{acknowledgments}

\section*{Codes}
The code used to generate the results for this paper will become available as a public repository at \url{https://gitlab.com/nicolasevangelou/ml_phase_transition} upon acceptance of the paper.

\appendix

\section{}\label{sec:Methods}
\subsection{Data Collection}
\label{sec:data_collection}
Given the ABM model described in Section \ref{sec:Model} in the main text, we generated trajectories at 18 equidistant values of $\sigma \in [0.5,2.2]$ \cite{zagli2023dimension}. 
For each value of $\sigma$ we sampled 100 trajectories for different initial conditions (agents' distribution) drawn from the Pearson distribution (implemented in Matlab) with prescribed mean, standard deviation, skewness, and kurtosis. The values for the mean, standard deviation, skewness, and kurtosis were chosen randomly from a prescribed grid of equidistant points: 
The mean was chosen from the range $[-2.0, 2.0]$, with increments of $0.2$. The standard deviation was selected from the range $[0.0,2.0]$, with increments of $0.1$. The skewness was chosen from the range $[-2.0,2.0]$, with increments of $0.2$ and the kurtosis was selected from the range $[0.0,15]$, with increments of $0.72$. This scheme ensured a dense sampling both, in parameter and state space.

The integration time of the ABM, containing $N=12,000$ agents, was set to $t_f=10$ with a time-step $dt=0.005$. Data were collected every five snapshots which led to a total of $400$ snapshots per trajectory. Therefore,  for a single value of the parameter, the total number of initial data is $n=40,000$ ($400 \times 100$) while for multiple values of the parameter $n=720,000$ ($400 \times 100 \times 18$). Note that a number of trajectories explodes and thus, are omitted from any further computations.

An additional preprocessing step was applied to the data before the Diffusion Maps computation that includes discarding a short transient $t_{cut}$ for each trajectory. This ensures that the fast transients have decayed and the collected data contains only the long-term dynamics (that live on the slow manifold). When dealing with multiple parameter values we chose $t_{cut}=1$ while for a single value of the parameter $\sigma=1$ we choose $t_{cut} = 0.5$ to make sure that enough transients are close to the unstable steady state, otherwise, the manifold would appear as two clusters.
As a test set, we sampled nine trajectories for each of the seven values of $\sigma = \{0.57, 0.85, 1.11, 1.75, 1.9, 2.06, 2.25 \}$ not included in the training set. The selection of these test values was made to validate the predictions of our Euler neural network for dynamics before and after the bifurcation. 

\subsection{Diffusion Maps}\label{sec:Dmaps}
The Diffusion Maps algorithm, introduced by Coifman and Lafon  \cite{coifman2006diffusion} can be used to discover a low-dimensional parameterization of high-dimensional data $\mathbf{X} = \{\vect{x}_i\}_{i=1}^n$ with each $\vect{x}_i \in \mathbb{R}^m$. Diffusion Maps constructs a weighted graph $\mathbf{K} \in \mathbb{R}^{n \times n}$ between the sampled data points by using a kernel function. A common choice, also used in our case, is the Gaussian kernel
\begin{equation}
    \label{eq:affinity_kernel}
    K(\vect{x}_i,\vect{x}_j) = \exp\bigg(\frac{- || \vect{x}_i - \vect{x}_j||_2^2}{2\varepsilon} \bigg),
\end{equation}
where $\varepsilon$ is a positive hyperparameter that controls the rate of the kernel's decay. The metric $|| \cdot ||$ in our case was chosen as the $\ell^{2}$ norm but different metrics are also possible.

To discover a low-dimensional manifold, regardless of the sampling density, the following normalization is required 
    \begin{equation}
    \mathbf{\Tilde{K}} = \mathbf{P}^{-1}\mathbf{K}\mathbf{P}^{-1}
\end{equation}

\noindent
where $P_{ii} = \sum_{j=1}^n K_{jj}$. A second normalization of $\mathbf{\Tilde{K}}$ recovers a row-stochastic, Markovian matrix

\begin{equation}
    \mathbf{M} = \mathbf{D}^{-1}\mathbf{\Tilde{K}}
\end{equation}

\noindent
where $\mathbf{D}$ is a diagonal matrix defined as  $D_{ii} = \sum_{j=1}^n \Tilde{K}_{ij}$.

\noindent
The entries of matrix $\mathbf{M}$ can be seen as probabilities of jumping from one point to the other. The eigendecomposition of $\mathbf{M}$,

\begin{equation}
    \mathbf{M}\vect{\phi}_i = \lambda_i\vect{\phi}_i,
\end{equation}

\noindent
provides a set of eigenvectors $\vect{\phi}_i$ and corresponding eigenvalues $\lambda_i$. To obtain a more parsimonious representation of the original data set $\mathbf{X}$ proper selection of the eigenvectors is needed. 
If the \textit{intrinsic} dimension of the data is small, this selection can be achieved by visual inspection of the non-harmonic eigenvectors (eigenvectors that span independent directions) \cite{evangelou2023double}. 
Alternatively, the local-linear regression algorithm proposed by Dsilva et al. \cite{dsilva2018parsimonious} can be used for selecting the non-harmonic eigenvectors. If, the number of independent non-harmonic eigenvectors is smaller than the dimension $m$ of the data $\mathbf{X}$, then dimensionality reduction has been achieved. In our work, the Python library \textit{datafold} \cite{lehmberg2020datafold} was used for the Diffusion Maps and the local-linear regression algorithms.

\subsubsection{Nyst\"om Extension}
\label{sec:nystrom}
The Nystr\"om Extension formula provides a numerical approximation of eigenfunctions  of the form \cite{fowlkes2004spectral}

\begin{equation}
    \int_{a}^b M(\vect{x}_j,\vect{x}_i)\phi_i(\vect{x}_i)d\vect{x}_i = \lambda \phi({\vect{x}_j}).
\end{equation}

In our work, Nystr\"om extension is utilized when new \textit{out-of-sample} data points are given, e.g., $\vect{x}_{new} \notin \mathbf{X}$. To generate the Diffusion Maps coordinates $\vect{\phi}_{new}$ Nystr\"om extension uses an interpolation scheme based on the kernel computations and normalizations applied during the dimensionality reduction step, discussed in the previous Section. The Nystr\"om extension formula reads

\begin{equation}
    \phi_i(\vect{x}_{new}) = \frac{1}{\lambda_i} \sum_{j=1}^n \tilde{\mathbf{M}}(\vect{x}_{new},\vect{x}_{j})\phi_i(\vect{x}_j),
\end{equation}

\noindent
where $\phi_i(x_{new})$ denotes the estimated $i$-th eigenvector for the data point $\vect{\phi}_{new}$, $\lambda_i$ denotes the corresponding eigenvalue, $\phi_i(\vect{x}_{j})$ denotes the $j$-th component of the $i$-th eigenvector and $\tilde{\mathbf{M}}(\cdot,\cdot)$ denotes the kernel function used to determine the similarity of $\vect{x}_{new}$ to all the points in $\mathbf{X}$.

\subsubsection{Diffusion Maps on ABM data} \label{sec:dmaps_ABM}
In this section, we provide details on how the Diffusion Maps algorithm was computed on histograms and moments from the ABM.  In the first case, for each snapshot of the ABM we constructed a histogram as an approximation of the agents' density. Each histogram contains $40$ equidistant bins defined in the range $[-4,4]$. This range ensures that all agents in the collected (training) data lie in between. Note,  that the method is insensitive to the selected number of bins.  To reduce the computational cost of the Diffusion Maps, we subsampled the training data uniformly \cite{lehmberg2020datafold} which resulted in about $3,500$ data points when Diffusion Maps applied for a single parameter and about $19,000$ data points for multiple values of the parameter. The hyperparameter $\varepsilon$ was selected as the square of the median of the pairwise distances multiplied by a constant $c$, where $c = 0.03$ for a single value of the parameter $\sigma$ and $c=20$ for multiple values of the parameter $\sigma$.  In the second case, Diffusion Maps was computed on the sampled moments, $M_1,M_2,M_3,M_4$.
Again, to alleviate the computational cost we subsampled the data \cite{lehmberg2020datafold},  which results in $N \sim 11,000$. The hyperparameter $\varepsilon$ was also selected for this case by computing the median of the pairwise distances multiplied with $c=10$.

\subsubsection{Diffusion Maps: Additional Results}
\label{sec:dmaps_SI}

\begin{figure}[ht!]
    \centering
\subfigure[]{\includegraphics[clip,width=0.18\textwidth]{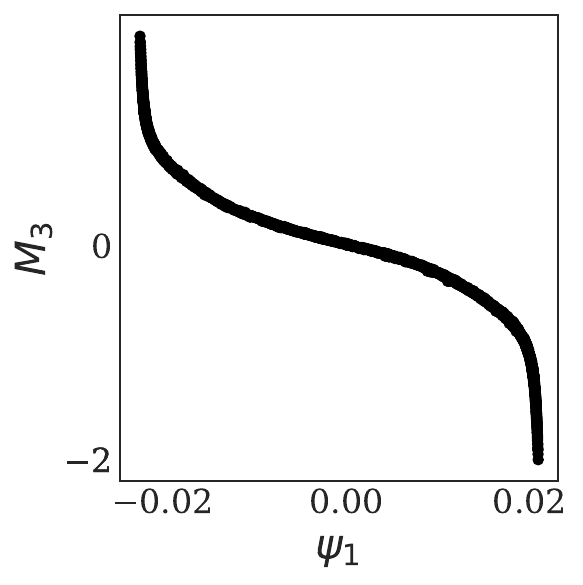}\label{fig:psi1_vs_m3}}
\subfigure[]{\includegraphics[clip,width=0.18 \textwidth]{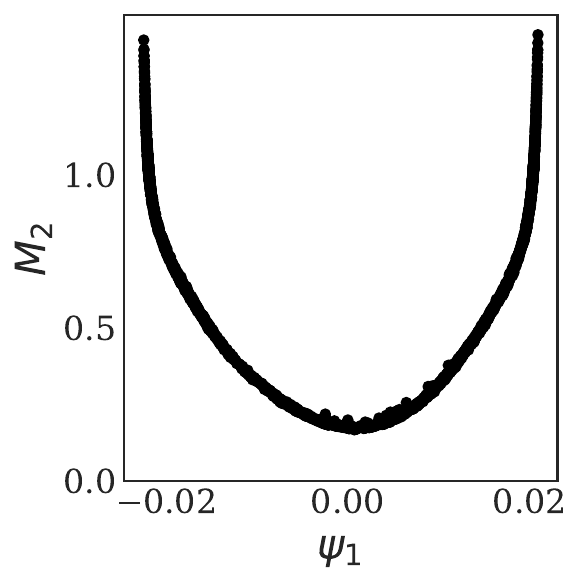}
\label{fig:psi1_vs_m2}}
\subfigure[]{\includegraphics[clip,width=0.18 \textwidth]{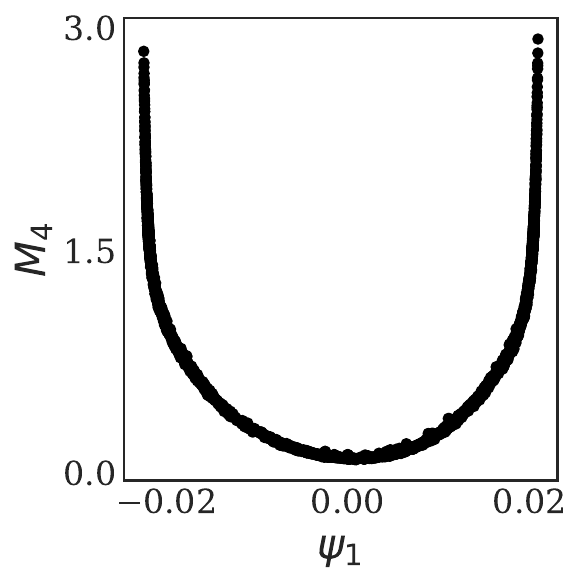}\label{fig:psi1_vs_m4}}

\subfigure[]{
{\includegraphics[clip,width=0.45\textwidth]{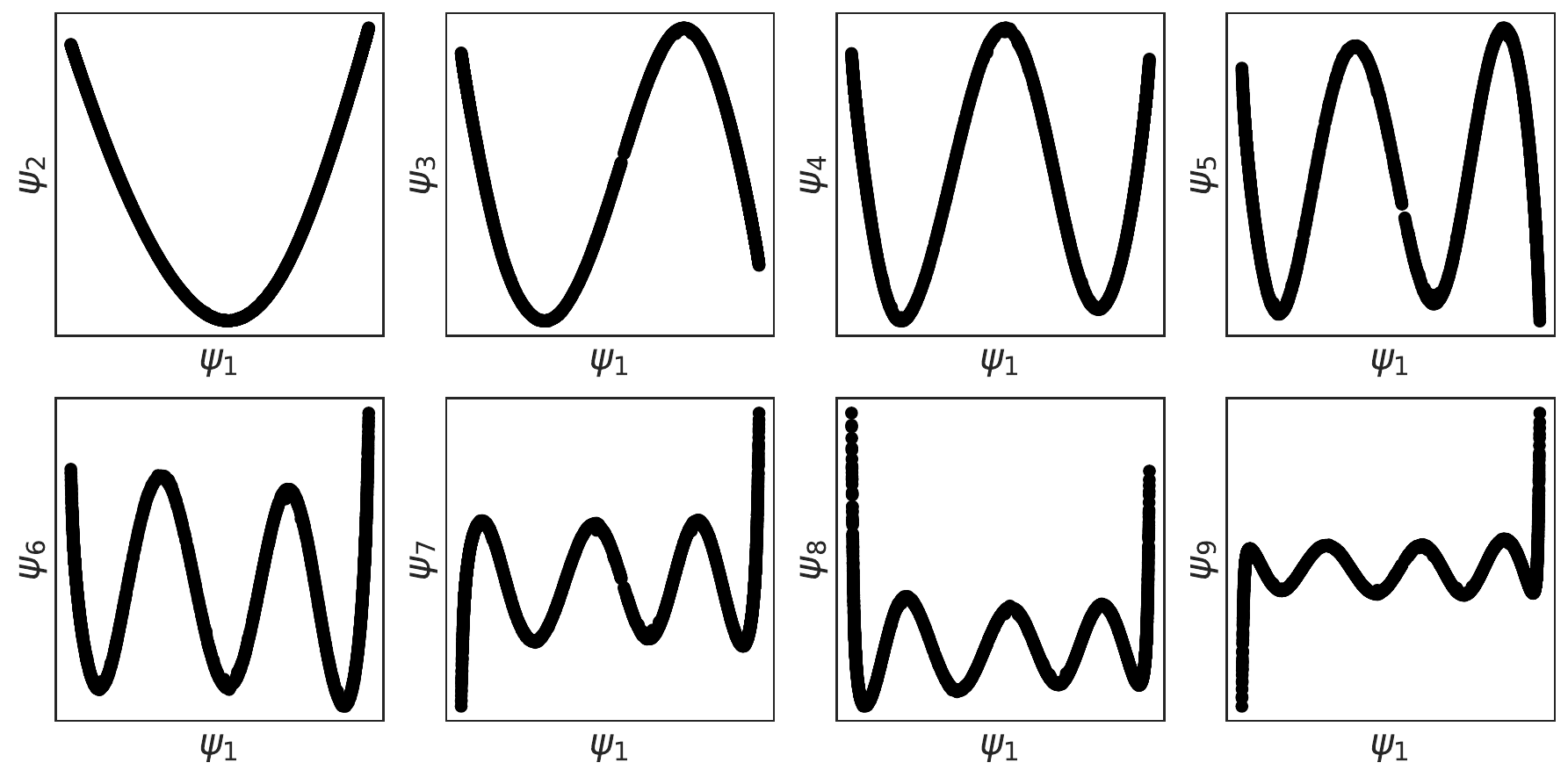}\label{fig:evec1_vs_evecs}}}

    \caption{(a-c) The Diffusion Maps coordinate $\psi_1$ is plotted against the computed moments $M_2$, $M_3$ and $M_4$ respectively.
    (d) The eigenvector $\psi_1$ is plotted against the eigenvectors $\psi_2- \psi_9$. This supports our argument that the eigenvectors $\psi_2- \psi_9$ are harmonics of $\psi_1$.}
    \label{fig:dmaps_1_parameter}
\end{figure}

In this section, we provide additional results for the Diffusion Maps computation for $\sigma=1$ on the histograms and the Diffusion Maps computations on the moments $M_1,M_2,M_3,M_4$. We showed in the main text that the Diffusion Maps coordinate $\psi_1$ is one-to-one with the first moment $M_1$. Here we show that $M_3$ is also one-to-one with $\psi_1$ (Figure \ref{fig:psi1_vs_m3}) and that the even moments $M_2$ and $M_4$ can be seen as functions of $\psi_1$ (Figures \ref{fig:psi1_vs_m2},\ref{fig:psi1_vs_m4}). 


To strengthen, the argument that the manifold in this case is one-dimensional we show in Figure \ref{fig:evec1_vs_evecs} the eigenvectors $\psi_2-\psi_9$ plotted against the first non-trivial eigenvector $\psi_1$. From Figure \ref{fig:evec1_vs_evecs} it appears that the eigenvectors $\psi_2-\psi_9$ are harmonics of $\psi_1$. This suggests that the manifold for a fixed value of $\sigma$ is one-dimensional. 

\begin{figure}[ht!]
    \centering
    \subfigure[]{ \includegraphics[width=0.14\textwidth]{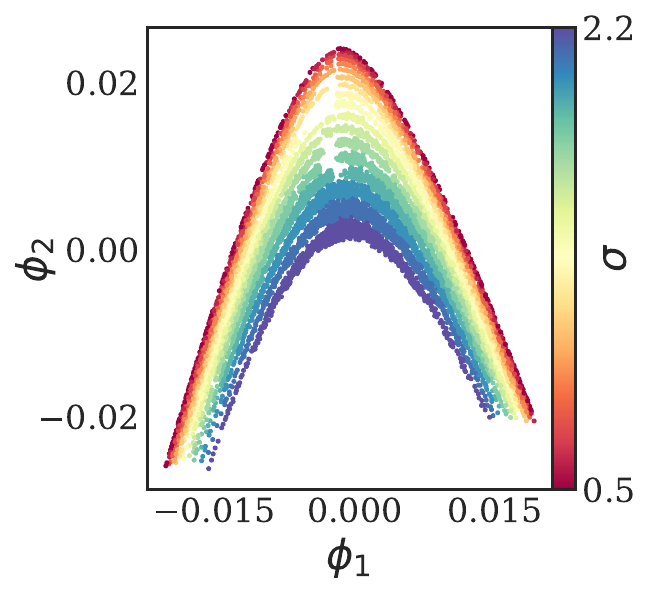}
\label{fig:non_harmonic_color_sigma_mom}}
 \subfigure[]{ 
\includegraphics[width=0.14\textwidth]{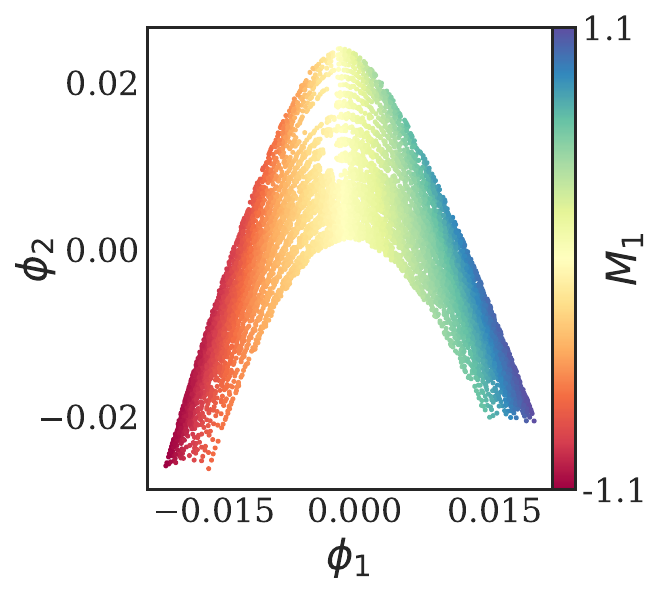}
\label{fig:non_harmominc_m1_mom}}
\subfigure[]{
\includegraphics[width=0.13\textwidth]{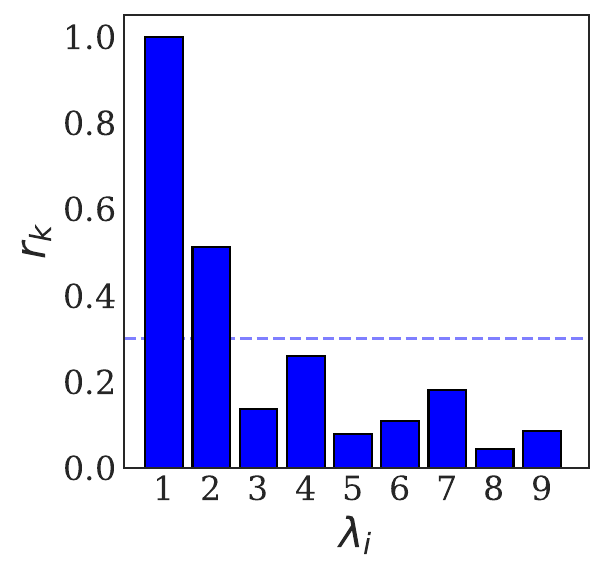}
\label{fig:local_linear_moments}}
    \caption{
    (a-b) The non-harmonic coordinates $\phi_1$ and $\phi_2$ colored with $\sigma$ and $M_1$ respectively. (c) The residual $r_k$ indicates $\phi_1$ and $\phi_2$ are the two non-harmonic coordinates.}
\label{fig:data_diffusion_maps_moments}
\end{figure}

The Diffusion Maps algorithm applied using the computed first four moments reveals a two-dimensional manifold embedded in four dimensions. The first two eigenvectors in this case $\phi_1$ and $\phi_2$, shown in Figures \ref{fig:non_harmonic_color_sigma_mom}
and \ref{fig:non_harmominc_m1_mom} are the non-harmonic, eigenvectors that parametrize the data. This is corroborated by the larger residuals $r_k$ of the first two eigenvectors shown in Figure \ref{fig:local_linear_moments}. Also in this case, the parameter $\sigma$ and the first moment $M_1$ appear visually as functions of the Diffusion Maps coordinates (Figures \ref{fig:non_harmonic_color_sigma_mom}
and \ref{fig:non_harmominc_m1_mom}).This pair of coordinates, obtained by Diffusion Maps on moments, could have been used for the CAE's training instead of the Diffusion Maps coordinates ($\psi_1,\psi_2$)obtained from the computations on the histograms; but we omit this for brevity.






\subsection{Y-shaped conformal autoencoder}
\label{sec:Y_shaped_CAE_details}
The Y-shaped conformal autoencoder was initially presented in \cite{evangelou2022parameter}. In this work, the Y-shaped conformal autoencoder consists of the three connected (sub)networks

\begin{align}
    & \text{Encoder : } (\psi_1,\psi_2) \mapsto (\nu_1, \nu_2) \\
    & \text{Decoder : } (\nu_1, \nu_2) \mapsto (\hat{\psi}_1, \hat{\psi}_2) \\
    & \text{Estimator : } \nu_1 \mapsto \hat{\sigma}.
\end{align}

\noindent
The Encoder receives as inputs the two Diffusion Maps coordinates $\psi_1$, $\psi_2$ and maps them to the latent variables $\nu_1,\nu_2$. The Decoder aims to reconstruct the Diffusion Maps coordinates from the latent variables $\nu_1$, $\nu_2$. The Estimator has as input the latent coordinate $\nu_1$ and aims to learn a map from $\nu_1$ to the parameter $\sigma$.

    
The loss function used to train the Y-shaped conformal autoencoder consists of three parts: (a) The loss function of the Encoder-Decoder (autoencoder) $\mathcal{L}_\text{ae}$ that aims to reconstruct the input itself  (b) the loss function of the Estimator, $\mathcal{L}_\text{est.}$, that aims to reproduce the parameter $\sigma$ given $\nu_1$ and (c) the loss function for imposing the conformality constrain, $\mathcal{L}_\text{con.}$, between $\nu_1$ and $\nu_2$,
\begin{equation}
\label{eq:conformality_1}
       \langle \nabla \nu_1   ,  \nabla \nu_2 \rangle = 0     
\end{equation}

\noindent
where the gradient $\nabla$ is in terms of the Diffusion Maps coordinates ($\psi_1$,$\psi_2$) of the input and $\langle \cdot, \cdot \rangle$ denotes the inner product between the two vectors. The gradients were computed by using the automatic differentiation of Pytorch \cite{pytorch}. In practice, instead of using equation \ref{eq:conformality_1}, one can minimize the angle between the vectors 

\begin{equation}
\label{eq:cosine}
    \text{cos}\theta = \frac{ \nabla\nu_1 \cdot \nabla \nu_2 }{||\nabla \nu_1|| ||\nabla \nu_2 || }
\end{equation}

\noindent
which stabilizes the training of the network. We describe details for training the network, specifics about the architecture used, and the choice of hyperparameters in the next section.

\subsubsection{Hyperparameter selection and training procedure}
The implementation of the Y-shaped conformal autoencoder was done with the Pytorch Python library \cite{pytorch}.

Each (sub)network (Encoder, Decoder, Estimator) in the architecture of the Y-shaped conformal autoencoder consists of five fully connected layers. The first four hidden layers have 20 neurons and tanh(t) activation functions and the fifth has no activation function (linear activation) and its size depends on the size of the desired output.  The ADAM optimizer was chosen for training the overall network. We chose mini-batches of size $32$ to train the network. The learning rate was selected as $\eta = 0.001$ and $500$ epochs. Training the network for a larger number of epochs leads to overfitting.

To train the network and test its generalization capability we used  about $19,000$ data points. We split the data into train$|$test$|$validation as $80$:$10$:$10$. We then rescaled the training data by using the \textit{MinMaxScaler} Python preprocessing scheme from \textit{sklearn}. We applied the same transformation for the validation and test sets. During the training of the network, we used only the training set to perform backpropagation and the validation set to get insight into the model's performance. The network did not \textit{see} the test set during training.

The optimization process we performed was \textit{heuristic}: for a fixed mini-batch two updates (backpropagation steps) were performed. The first step updates the weights of the Encoder-Decoder and the second step the weights of the Estimator-Encoder. Altering this training protocol is possible. In Algorithm \ref{alg:conformal-autoencoder} below we provide a more detailed description of the network's training.

Upon training of the neural network, the estimated MSE for the autoencoder's reconstruction, was $\mathcal{L}_{\text{ae}} = 1.36e-08$ on the train set and $\mathcal{L}_{\text{ae}} = 1.37e-08$ on the test set. The MSE for the Estimator was $\mathcal{L}_{\text{est.}} = 1.2e-03$ for the train set and $\mathcal{L}_{\text{est.}} =1.2e-03$ for the test. The average value of the $\text{cos}\theta$, on the train set was $\mathcal{L}_\text{con.}=2.72e-05$ and on the test set $\mathcal{L}_\text{con.}=2.63e-05$.

\begin{algorithm}[htp!]
    \caption{The algorithm illustrates a full iteration during training of the Y-shaped conformal autoencoder. We set the scale parameter to $\alpha=10$. The learning rate is denoted as $\eta$.}
    \textbf{\underline{Input}}: Diffusion Maps coordinates $\psi_1$,$\psi_2$ and parameter values $\sigma$.\\
    \textbf{\underline{Output}:} The weights of: (i) Encoder ($\theta_{\text{encoder}}$), \\ (ii) Decoder ($\theta_{\text{decoder}}$), (iii) Estimator ($\theta_{\text{estimator}}$).\\
    For {$i=1,2,...,T$}
        \begin{enumerate}
             \item Predict:
             \begin{align*}
                 (\nu_1,\nu_2) &=\text{Encoder}(\psi_1,\psi_2) \\
             (\hat{\psi}_1,\hat{\psi}_2)&=\text{Decoder}(\nu_1,\nu_2)
             \end{align*}
            \item Compute Autoencoder (Encoder-Decoder) and Conformality Losses: 
            \begin{align*}
                \mathcal{L}_1 &= \mathcal{L}_{\text{ae}} +  \mathcal{L}_{\text{con.}} \\
                &= \text{MSE}(\hat{\vect{\psi}},\vect{\psi}) + \alpha \text{MSE}(\text{cos}\theta,0)
            \end{align*}
            \item Backpropagation step - update weights (illustration with gradient descent):
            \begin{align*}
                &\theta_{\text{encoder}} -= \eta{\theta_{\text{encoder}}} \mathcal{L}_1 \\
                & \theta_{\text{decoder}} -= \eta{\theta_{\text{decoder}}} \mathcal{L}_1
            \end{align*}
             \item  Predict:
             \begin{align*}
                &(\nu_1,\nu_2) = \text{Encoder}(\psi_1,\psi_2) \\
                &\hat{\sigma}  = \text{Estimator}(\nu_1)
             \end{align*}
             \item Compute Estimator Loss:
             \begin{align*}
                 \mathcal{L}_{\text{est.}}  = \text{MSE}(\hat{\sigma}, \sigma)
             \end{align*}
             \item Backpropagation step - update weights (illustration with gradient descent)
            \begin{align*}
            &\theta_{\text{estimator}} -= \eta{\theta_{\text{estimator}}} \mathcal{L}_{\text{est.}} \\
            &\theta_{\text{encoder}} -= 
\eta{\theta_{\text{encoder}}} \mathcal{L}_{\text{est.}}
            \end{align*}
        \end{enumerate}
    \label{alg:conformal-autoencoder}
    \hrulefill
\end{algorithm}

\subsubsection{Y-shaped conformal autoencoder: Additional Results}
\label{sec:Y-shaped_SI}

\begin{figure}[ht!]
    \centering 
    \subfigure[]{
\includegraphics[width=0.2\textwidth]{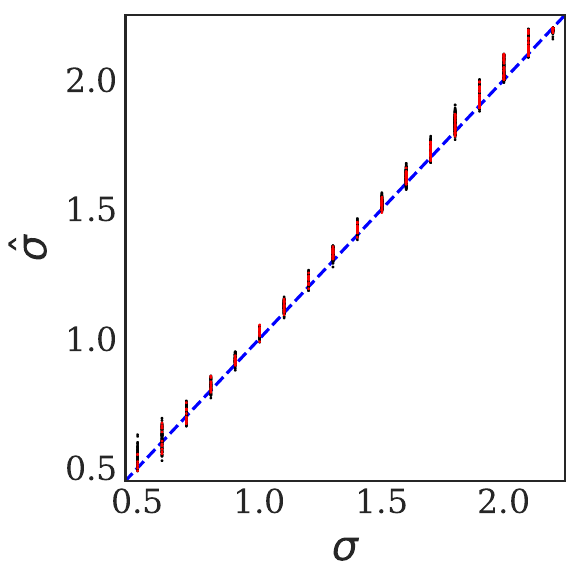
} \label{fig:cae_predict_sigma}}

\subfigure[]{
\includegraphics[width=0.21\textwidth]{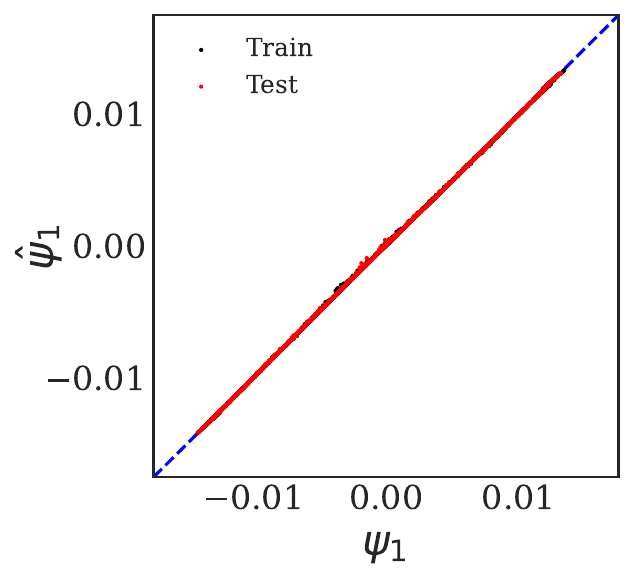}\label{fig:cae_predict_psi1}}
\subfigure[]{
\includegraphics[width=0.21\textwidth]{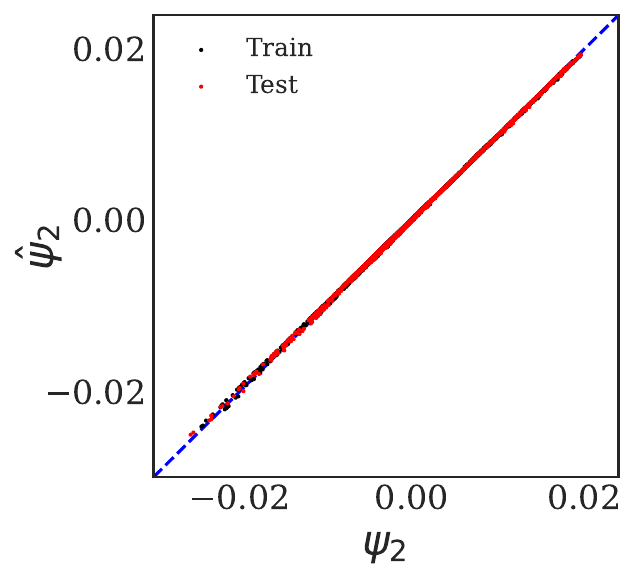}\label{fig:cae_predict_psi2}}
    \caption{
    (a) The true values of the parameter $\sigma$ are plotted against the reconstructed by the Estimator for train (black) and test (red) points. The blue dashed line indicates $y=x$. (b-c) The true values of the Diffusion Maps coordinates $\psi_1$ and $\psi_2$ are plotted against the reconstructed $\hat{\psi}_1$, $\hat{\psi}_2$ coordinates by the autoencoder. }
\end{figure}

In this section, we provide additional results regarding the Y-shaped conformal autoencoder described in the main text. 
%
The ability of the Estimator to predict the parameter $\sigma$ from $\nu_1$ is shown in Figure \ref{fig:cae_predict_sigma} for train (black dots) and test (red dots) points.
We also illustrate, in Figures \ref{fig:cae_predict_psi1}, \ref{fig:cae_predict_psi2}, the reconstruction of the autoencoder 
for train and test points.

\subsection{Forward Euler neural network}
\label{sec:euler_nn_methodology}
In this section, we describe how we identified the right-hand side of an ordinary differential equation directly from data. Let $\nu_2(t)$ be a state variable whose dynamics are governed by a $\sigma$-dependent ODE given by the general form
\begin{equation}
\label{eq:ode_equation}
    \dot{\nu}_2(t) = f(\nu_2(t);\sigma)
\end{equation}
 Our goal is to construct a neural network architecture inspired by numerical integrators of ODEs to estimate the right-hand-side $f(\nu_2(t);\sigma)$. To this end, we constructed a forward Euler residual neural network depicted in Figure \ref{fig:euler_nn}. To train this network, we do not require long trajectories but only snapshots of the form $\mathcal{D} = \{\nu_2(t+h), \nu_2(t), \sigma, h \}$ where $\nu_2(t)$ is the state variable at time $t$, $\nu_2(t+h)$ is the state variable after a small time-step $h$. Given sampled data in the form of $\mathcal{D}$ we wish to approximate $f$ by using a neural network, with weights denoted as  $\theta$. 

To formulate the loss used to train the network $f_\theta$ we remind the reader that the forward Euler approximates the evolution of an ODE, by a small positive step $h$, as 
\begin{equation}
    \nu_2(t+h) = \nu_2(t) + hf(\nu_2;\sigma).
\end{equation}
In our case, without having access to $f$ but only data in the form $\mathcal{D}$ we wish to approximate the right-hand-side with the neural network $f_\theta$. This is achieved by performing for each pair of inputs ($\nu_2(t),\sigma$) one integration step of size $h$, estimating the evolved dynamics \hbox{$\hat{\nu_2}(t + h)$} and minimizing the loss
\begin{equation}
    \label{eq:loss_nn_euler}
    \mathcal{L}(\theta| \nu_2(t), \nu_2(t+h), h, \sigma) = || \hat{\nu}_2(t+h) - \nu_2(t+h) ||^2.  
\end{equation}
For our computations, the time step $h$ was kept constant but the overall approach can be easily extended to handle also varying time steps $h$.

\subsubsection{Hyperparameter selection and training procedure}
The implementation of the forward Euler neural networks was done with the Tensorflow/Keras Python libraries \cite{tensorflow2015}.

To train the forward Euler neural network we need to ensure our data is in the form of snapshots $\mathcal{D} $. To achieve that we used the Nystr\"om extension formula (see Section \ref{sec:nystrom}) to all the available sampled trajectories and obtained the corresponding trajectories in $\psi_1\psi_2$. We then evaluated the Encoder to get trajectories in terms of $\nu_2$. These two steps provided us with about $600,000$ snapshots. We then split the data into train$|$test$|$validation as 80$:$10$:$10. We then centered and whitened the data (based on the mean and variance of the training set) and applied the same transformation to the validation and test set.

The architecture consisted of two hidden layers with $10$ neurons each. The first hidden layer had a tanh(t) activation function and the second hidden layer had a linear activation function. We used ADAM to optimize this network. The mini-batch size was set to $32$, the number of total epochs to $100$, and the learning rate to $\eta = 0.001$. The learning curves for the training and validation are shown in Figure \ref{fig:loss_function_euler_nn}. The ability of the network to fit the training set and generalize is shown in \ref{fig:euler_nn_predictions}. 
Upon training of the network, the MSE on the train and test sets were $2.63e-03$ and $2.79e-05$ respectively. This neural network was used for the computations reported in \ref{sec:Euler_NN_one}. We did not record the MSE for all the $5000$ networks used to check the robustness of our approach in identifying the critical transition.

\subsubsection{Euler neural network: Additional Results}
\begin{figure}[ht!]
    \centering
\subfigure[]{\includegraphics[width=0.21 \textwidth]{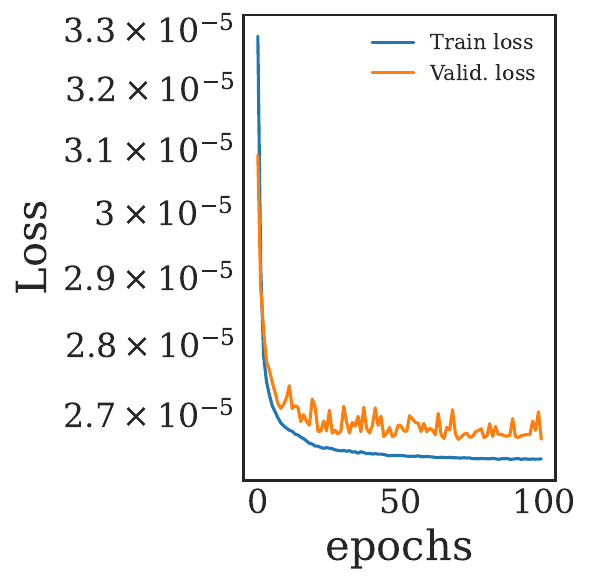} \label{fig:loss_function_euler_nn}}
    \subfigure[]{
\includegraphics[width=0.215\textwidth]{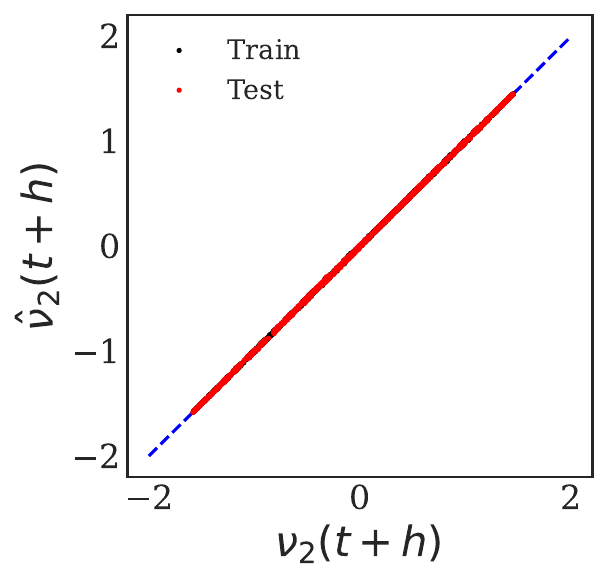} \label{fig:euler_nn_predictions}}
    \caption{(a) The learning curves for the train and validation sets. (b) The true values of the state variable $\nu_2(t+h)$ plotted against the predicted by the Euler neural network $\hat{\nu}_2(t+h)$ for the train (black points) and test (red points). }
\end{figure}

In this section, we provide some additional results for the forward Euler neural network.
 The learning curves for the training and validation are shown in Figure \ref{fig:loss_function_euler_nn}. The ability of the network to fit the train set and generalize is shown in Figure \ref{fig:euler_nn_predictions}.
 Note that the normalized value of the state variable $\nu_2$ is shown in Figure \ref{fig:euler_nn} and that the test points in this include values of the parameter $\sigma$ that are also in the training set.
 
\subsection{Separation of time scales}
\label{sec:sep_time_scales}
\begin{figure}[ht!]
    \centering
{\includegraphics[clip,width=0.35 \textwidth]{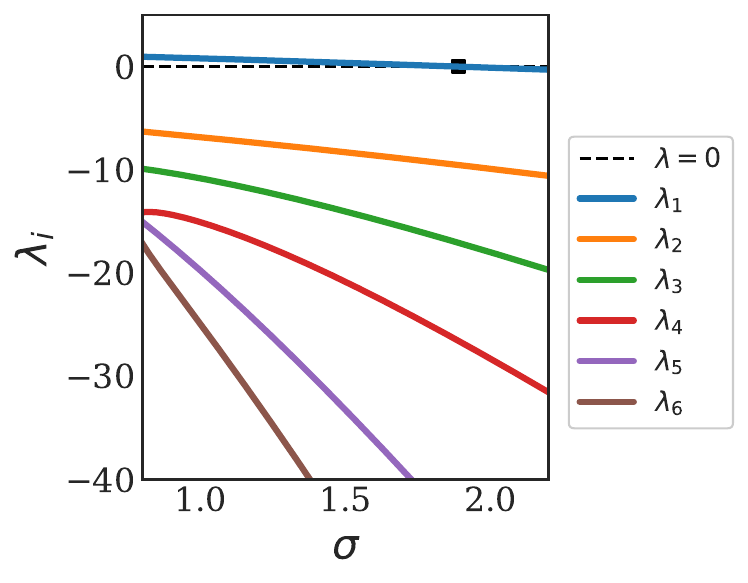}}
    \caption{The eigenvalues of the Jacobian computed across a range of parameter values $\sigma$ for the moments' equations with six ODEs presented in \cite{zagli2023dimension}. The slowest eigenvalue (upper line) $\lambda_1$ crosses the real axis at $\sigma=1.89$ where the bifurcation occurs. The black square indicates the value of $\sigma$ where the (slowest) eigenvalue crosses zero.}
    \label{fig:time-scale-separation}
\end{figure}

In this Section, we discuss the separation of time scales for the reduced order models based on the equations for the moments proposed in \cite{zagli2023dimension}. As shown in Figure \ref{fig:time-scale-separation} the slowest eigenvalue $\lambda_1$ is at least one order of magnitude smaller than the second slowest $\lambda_2$. This suggests that the dynamics after a short transient are slaved in the direction of this slowest eigenvalue and are effectively one-dimensional.

\clearpage
\bibliography{apssamp}

\end{document}